# On series solutions of Volterra equations

S. A. Belbas
Mathematics Department
University of Alabama
Tuscaloosa, AL. 35487-0350.

e-mail: SBELBAS@GP.AS.UA.EDU

Abstract.  We derive formulae for the calculation of Taylor coefficients of solutions to systems of Volterra integral equations, both linear and nonlinear, either without singularities or with singularities of Abel type and logarithmic type. We also obtain solutions to certain systems of Volterra equations of the first kind. In all cases except the case of logarithmic singularities, we obtain recursive formulae for the calculation of the Taylor coefficients. In certain cases, we give proofs of convergence and rigorous estimates of the radius of convergence.

1. Introduction.

In recent years, there has been increasing interest in Taylor series solutions of Volterra and Fredholm integral equations and integrodifferential equations. We mention the papers [B, RZQ, W, Y, YS].

Taylor methods have the useful feature that, whenever they are applicable, they produce the exact Taylor coefficients of the solution; this is an advantage in case an approximation has been obtained and later higher accuracy is wanted: it is not necessary, in that case, to solve the problem anew, but the coefficients that were initially calculated can be used as input in the calculation of additional coefficients; it is also an advantage in the sense that the Taylor coefficients are intrinsic characteristics of the solution, and do not depend on any choices made in the implementation of an algorithm (by contrast, approximations based on numerical integration, finite differences, and other discretizations, depend on the choice of the discretization grid and on the nature of the discretization). Furthermore, it is always useful, for comparison and validation of algorithms, to have as many different ways as possible for the approximate solution of systems of Volterra integral equations.

In the past, series expansion methods did not receive a lot of attention as methods for finding approximate solutions to integral equations, due to the fact that such methods require the calculation of derivatives, which used to be an undesirable feature for numerical methods. The well-known recent developments in software for automatic differentiation has led to reconsideration of series expansion methods for integral equations, and consequently several papers have appeared on this subject. The relevance of automatic differentiation software to series expansion methods for integral equations has already been noted in [G].

It is natural to ask, whether a method similar to the well-known techniques for ordinary differential equations (classic methods associated with the names of Frobenius, Fuchs, Schlesinger, Floquet, et al.) can also be used for systems of Volterra integral equations, with or without singularities, linear or nonlinear. The ODE techniques rely on 3 main ingredients: (i) differentiation of functional series; (ii) multiplication of series; (iii) choice of a system of basis functions that remains closed under the operations of differentiation and multiplication of series, to the extent that these two operations are relevant to solving a differential equation. We can ask whether an analogous method, with "integration" substituted in place of "differentiation" in feature (i) above, can be applied to systems of Volterra integral equations. We shall show, in this paper, that the answer is affirmative. A fourth feature of series expansion methods for linear differential equations is the fact that they lead to recursive relationships for the determination of the unknown coefficients of a solution; this property sometimes, but not always, extends to systems of Volterra integral equations.

In this paper, we provide the recursive or, in some cases, non-recursive formulae for the series solution of regular systems of Volterra equations of the second kind, certain cases of Volterra systems of the first kind, systems of Volterra equations with Abel-type singularities, and systems of Volterra equations with logarithmic singularities. Whenever possible, we also provide a rigorous analysis of the convergence and the radius of convergence of the relevant functional series.

Technical note: We use two notational conventions. First, a summation over an empty set of indices will be, by definition, equal to 0; this convention allows us to write equations without explicitly dealing with special cases when a condition defining a set of indices becomes impossible. Second, when using subscripts or superscripts, we may separate indices by commas (when the indices become algebraic expressions), or simply write indices next to each other without commas (when the indices are single letters); for example, the indices "n" and "n–j" in "$K_{n,n-j}$" have the same status as the indices "α" and "β" in "$K_{\alpha\beta}$".

## 2. Outline of the differentiation approach.

In this preliminary section, we outline, for the sake of completeness, an approach that is conceptually due to [G]. The problem we consider in this section is more general than the problem of [G], and the exposition that follows is different from that in [G].
We consider the nonlinear Volterra integral system

$$x(t) = \xi(t) + \int_0^t f(t,s,x(s))ds$$

--- (2.1)

with vector-valued unknown function x(t). We denote by N the dimension of x(t), i.e. x(t) takes values in $\mathbb{R}^N$. We shall outline an approach based on the method introduced in [G]; although this method is due to [G], our description is modified to fit the purposes of the present paper. The problem treated in [G] was an initial-value problem for an integro-differential equation with scalar-valued unknown function; it is clear that the vector-valued model (2.1) also includes the integro-differential equation of [G].
For the purposes of this section, we assume that all functions are sufficiently many times continuously differentiable, and we shall not specify the degree of differentiability that is needed in every calculation, since that information is contained in each formula that involves derivatives.
Assuming that x(t) is sufficiently many times differentiable near t=0, we want to find a Taylor polynomial approximation to x(t), i.e. we want to find the coefficients $X_n$ in a finite Taylor expansion with remainder,

$$x(t) = \sum_{n=0}^{m} t^n X_n + O(t^{m+1})$$

--- (2.2)

If x(t) is real analytic, the same method produces any finite number of terms of the Taylor series for x(t),

$$x(t) = \sum_{n=0}^{\infty} t^n X_n$$

--- (2.3)

In fact, the calculations below will show that, if f is sufficiently many times differentiable, then so is also x(t).

First, we need to introduce some notation.
We set

$$f_t^{(n)}(t,s,x(s)) := \frac{\partial^n f(t,s,x(s))}{\partial t^n}, \quad f_t^{(n)}(t,t,x(t)) := \frac{\partial^n f(t,s,x(s))}{\partial t^n}\bigg|_{s=t}$$

--- (2.4)

The symbol $D_t^m$ will denote the m-th order total derivative with respect to t; for instance,

$$D_t^1 f(t,t,x(t)) = \left[\frac{\partial f(t,s,x(s))}{\partial t} + \frac{\partial f(t,s,x(s))}{\partial s} + \nabla_x f(t,s,x(s)) \bullet \frac{dx(s)}{ds}\right]\bigg|_{s=t},$$

$$D_t^2 f(t,t,x(t)) = \left[\frac{\partial^2 f(t,s,x(s))}{\partial t^2} + \frac{\partial^2 f(t,s,x(s))}{\partial s^2} + \left(\frac{dx(s)}{ds}\right)^T (H_x f(t,s,x(s)))\left(\frac{dx(s)}{ds}\right) + \right.$$

$$+ 2\nabla_x \left(\frac{\partial f(t,s,x(s))}{\partial t} + \frac{\partial f(t,s,x(s))}{\partial s}\right) \bullet \frac{dx(s)}{ds} + 2\frac{\partial^2 f(t,s,x(s))}{\partial s \partial t} +$$

$$\left. + \nabla_x f(t,s,x(s)) \bullet \frac{d^2 x(s)}{ds^2}\right]\bigg|_{s=t}$$

and so on.
Now, successive differentiations of (2.1), assuming all functions involved are sufficiently many times differentiable, gives

$$\frac{dx(t)}{dt} = \frac{d\xi(t)}{dt} + f(t,t,x(t)) + \int_0^t f_t^{(1)}(t,s,x(s))\,ds;$$

$$\frac{d^2 x(t)}{dt^2} = \frac{d^2 \xi(t)}{dt^2} + D_t^1 f(t,t,x(t)) + f_t^{(1)}(t,t,x(t)) + \int_0^t f_t^{(2)}(t,s,x(s))\,ds;$$

....

....

$$\frac{d^n x(t)}{dt^n} = \frac{d^n \xi(t)}{dt^n} + \sum_{r=0}^{n-1} D_t^r f_t^{(n-r-1)}(t,t,x(t)) + \int_0^t f_t^{(n)}(t,s,x(s))\,ds$$

--- (2.5)

For every sufficiently many times differentiable functions $\psi = \psi(t)$ and $\varphi = \varphi(t, x(t))$, we introduce the following notation:

$$\mathbf{D}_t^m \psi = \{\frac{d^p \psi}{dt^p} : p \leq m\};$$

$$\mathbf{D}_{t,x}^m \varphi = \{\frac{\partial^{p_0 + p_1 + \ldots + p_N} \varphi}{\partial t^{p_0} \partial x_1^{p_1} \ldots \partial x_N^{p_N}} : p_0 + p_1 + \ldots + p_N \leq m\}$$

--- (2.6)

Now, every derivative $D_t^r f_t^{(m)}(t, t, x(t))$ can be expressed in the form

$$D_t^r f_t^{(m)}(t, t, x(t)) = F(t, r, \mathbf{D}_{t,x}^r f_t^{(m)}(t, t, x(t)), \mathbf{D}_t^r x(t))$$

--- (2.7)

for a suitable function F. Eq. (2.7) is intended to show the quantities that affect $D_t^r f_t^{(m)}(t, t, x(t))$, not to provide a rule for calculating $D_t^r f_t^{(m)}(t, t, x(t))$.
The validity of (2.7) can be shown by induction. The explicit form for the function F can be obtained via multivariable extensions of the Faà di Bruno theorem for higher-order derivatives of composite functions [M, CS, P].
In view of (2.7), we can write (2.5) as

$$\frac{d^n x(t)}{dt^n} = \frac{d^n \xi(t)}{dt^n} + \Phi(t, n-1, \{\mathbf{D}_{t,x}^r f_t^{(n-r-1)}(t, t, x(t)) : 0 \leq r \leq n-1\}, \mathbf{D}_t^{n-1} x(t)) +$$
$$+ \int_0^t f_t^{(n)}(t, s, x(s)) \, ds$$

--- (2.8)

where

$$\Phi(t, n-1, \{\mathbf{D}_{t,x}^r f_t^{(n-r-1)}(t, t, x(t)) : 0 \leq r \leq n-1\}, \mathbf{D}_t^{n-1} x(t)) =$$
$$= \sum_{r=0}^{n-1} F(t, r, \mathbf{D}_{t,x}^r f_t^{(n-r-1)}(t, t, x(t)), \mathbf{D}_t^r x(t))$$

--- (2.9)

By taking t=0 in (2.9), we obtain

$$\frac{d^n x(t)}{dt^n}\bigg|_{t=0} = \frac{d^n \xi(t)}{dt^n}\bigg|_{t=0} + \Phi(0, n-1, \{\mathbf{D}_{t,x}^r f_t^{(n-r-1)}(0,0,x(0)): 0 \le r \le n-1\}, \mathbf{D}_t^{n-1} x(0))$$

--- (2.10)

Thus we have a recursive equation for the calculation of the quantities $x^{(n)}(0)$, and then we can find $X_n$ by using $X_n = \dfrac{x^{(n)}(0)}{n!}$.

This method can be easily specialized to the case of linear Volterra systems, of the form

$$x(t) = \xi(t) + \int_0^t k(t,s) x(s) ds$$

--- (2.11)

In the linear case, $f(t, t, x(t)) = k(t,t) x(t)$, thus

$$f_t^{(m)}(t, t, x(t)) = k_t^{(m)}(t, t) x(t)$$

--- (2.12)

and consequently

$$D_t^r f_t^{(n-r-1)}(t, t, x(t)) = \sum_{j=0}^{r} \binom{r}{j} \left(D_t^{r-j} k_t^{(n-r-1)}(t, t)\right) x^{(j)}(t)$$

--- (2.13)

Thus (2.10) becomes, in the linear case,

$$x^{(n)}(0) = \xi^{(n)}(0) + \sum_{j=0}^{n-1} \sum_{r=j}^{n-1} \binom{r}{j} \left(D_t^{r-j} k_t^{(n-r-1)}(0,0)\right) x^{(j)}(0)$$

--- (2.14)

### 3. A method akin to the Fuchs - Frobenius method in the regular case.

We consider again the system (2.1). We need some notation and definitions.

We shall denote by boldface letters N - dimensional multi-indices, for example $\mathbf{k} = (k_1, k_2, ..., k_N)$. The relationship $\mathbf{k} \geq \mathbf{l}$ means that $k_i \geq l_i \ \forall i = 1, 2, ..., N$. The N - dimensional multi-index with all entries zero will be denoted by $\mathbf{0}$. If x is a vector in $\mathrm{IR}^N$ with components $x_{(i)}$, $i = 1, 2, ..., N$, the power $x^{\mathbf{k}}$ is defined by $x^{\mathbf{k}} := \prod_{i=1}^{N} (x_{(i)})^{k_i}$. We assume the representations

$$\xi(t) = \sum_{i=0}^{\infty} t^i \Xi_i \ ; \ f(t, s, x) = \sum_{i=0}^{\infty} \sum_{j=0}^{\infty} \sum_{\mathbf{k} \geq \mathbf{0}} t^i s^j x^{\mathbf{k}} F_{ij\mathbf{k}}$$

--- (3.1)

We seek a solution $x(t)$ in the form

$$x(t) = \sum_{i=0}^{\infty} t^i X_i$$

--- (3.2)

In (3.2), each $X_i$ is an N - dimensional vector; we denote the components of $X_i$ by $X_{i(j)}$, $j = 1, 2, ..., N$. Then we have

$$(x(t))^{\mathbf{k}} = \left( \sum_{i=0}^{\infty} t^i X_i \right)^{\mathbf{k}} = \prod_{j=1}^{N} \left( \sum_{i=0}^{\infty} t^i X_{i(j)} \right)^{k_j} =$$

$$= \sum_{i_1=0}^{\infty} \sum_{i_2=0}^{\infty} \cdots \sum_{i_N=0}^{\infty} t^{\left( \sum_{j=1}^{N} i_j k_j \right)} \prod_{j=1}^{N} (X_{i(j)})^{k_j} =$$

$$= \sum_{l=0}^{\infty} t^l \left( \sum_{(i_1, i_2, ..., i_N): \ i_1 k_1 + i_2 k_2 + ... + i_N k_N = l} \prod_{j=1}^{N} (X_{i_j(j)})^{k_j} \right)$$

--- (3.3)

We define

$$Y_l := \{X_i : 0 \leq i \leq l\};$$

$$Z_{\mathbf{k}}(Y_l) := \sum_{(i_1, i_2, \ldots, i_N): i_1 k_1 + i_2 k_2 + \ldots + i_N k_N = l} \prod_{j=1}^{N} \left(X_{i_j(j)}\right)^{k_j};$$

$$\mathbf{i} := (i_1, i_2, \ldots, i_N); \quad \mathbf{i} \bullet \mathbf{k} := i_1 k_1 + i_2 k_2 + \ldots + i_N k_N;$$

$$X_{\mathbf{i}} := (X_{i_1(1)}, X_{i_2(2)}, \ldots, X_{i_N(N)}); \quad (X_{\mathbf{i}})^{\mathbf{k}} := \prod_{j=1}^{N} \left(X_{i_j(j)}\right)^{k_j};$$

$$I_l := \{\mathbf{k} : \exists \mathbf{i} \text{ such that } \mathbf{i} \bullet \mathbf{k} = l\}$$

--- (3.4)

Clearly, the second equation in (3.3) can be written in compact form as

$$Z_{\mathbf{k}}(Y_l) = \sum_{\mathbf{i}: \mathbf{i} \bullet \mathbf{k} = l} (X_{\mathbf{i}})^{\mathbf{k}}$$

--- (3.5)

The expression $\sum_{\mathbf{i}: \mathbf{i} \bullet \mathbf{k} = l} (X_{\mathbf{i}})^{\mathbf{k}}$ is indeed a function of $Y_l$, because $(X_{\mathbf{i}})^{\mathbf{k}} \equiv \prod_{j=1}^{N} \left(X_{i_j(j)}\right)^{k_j}$ depends only on those $X_{i_j(j)}$ for which $k_j \geq 1$ (since, for $k_j = 0$, we have $\left(X_{i_j(j)}\right)^{k_j} = \left(X_{i_j(j)}\right)^0 = 1$), and, for $k_j \geq 1$, the condition $\mathbf{i} \bullet \mathbf{k} = l$ implies $i_j \leq l$.

Eq. (3.3) can now be written as

$$\left(\sum_{i=0}^{\infty} t^i X_i\right)^{\mathbf{k}} = \sum_{l=0}^{\infty} t^l Z_{\mathbf{k}}(Y_l)$$

--- (3.6)

Consequently,

$$f(t, s, x(s)) = \sum_{i=0}^{\infty} \sum_{j=0}^{\infty} \sum_{l=0}^{\infty} \sum_{k \geq 0} t^i s^{j+l} Z_k(Y_l) F_{ijk}$$

--- (3.7)

Substitution of the series expansions into (2.1) yields

$$\sum_{i=0}^{\infty} t^i X_i = \sum_{i=0}^{\infty} t^i \Xi_i + \sum_{i=0}^{\infty} \sum_{j=0}^{\infty} \sum_{l=0}^{\infty} \sum_{k \geq 0} \frac{t^{i+j+l+1}}{j+l+1} Z_k(Y_l) F_{ijk}$$

--- (3.8)

and, after some changes of the indices of summation and the order of summation, we obtain the recursive equations

$$X_n = \Xi_n + \sum_{l=0}^{n-1} \sum_{i=0}^{n-l-1} \sum_{k \in I_l} \frac{1}{n-i} Z_k(Y_l) F_{i,n-i-l-1,k}$$

--- (3.9)

(The set $I_l$ has been defined in (3.4).)

Because the Taylor coefficients of x(t) are uniquely defined, it follows that (3.9) produces the same results as (2.10), using the identification $X_n \equiv \frac{x^{(n)}(0)}{n!}$. Consequently, (3.9) also produces the correct coefficients of a Taylor polynomial approximation to x(t), if the functions f and ξ are sufficiently many times continuously differentiable in all their arguments but not necessarily analytic.

The method just presented can, of course, be specialized to the case of a linear Volterra system. We consider the system (2.11), and we assume

$$x(t) = \sum_{i=0}^{\infty} t^i X_i \,;\, \xi(t) = \sum_{i=0}^{\infty} t^i \Xi_i \,;\, k(t,s) = \sum_{i=0}^{\infty} \sum_{j=0}^{\infty} t^i s^j K_{ij}$$

--- (3.10)

Substitution of (3.10) into (2.11) leads to recursive equations for the determination of the unknown coefficients $X_i$. First, we get

$$\sum_{i=0}^{\infty} t^i X_i = \sum_{i=0}^{\infty} t^i \Xi_i + \sum_{i=0}^{\infty} \sum_{j=0}^{\infty} \sum_{l=0}^{\infty} \frac{t^{i+j+l+1}}{j+l+1} K_{ij} X_l$$

--- (3.11)

With a change in the indices of summation and a change in the order of summation in the triple series, (3.11) can be written as

$$\sum_{n=0}^{\infty} t^n X_n = \sum_{n=0}^{\infty} t^n \Xi_n + \sum_{n=0}^{\infty} t^n \left( \sum_{l=0}^{n-1} \sum_{i=0}^{n-l-1} \frac{1}{n-i} K_{i,n-i-l-1} X_l \right)$$

--- (3.12)

from which we get the recursive equation (discrete-time Volterra equation)

$$X_n = \Xi_n + \sum_{l=0}^{n-1} \sum_{i=0}^{n-l-1} \frac{1}{n-i} K_{i,n-i-l-1} X_l$$

--- (3.13)

It is seen from (3.13) that, if the Taylor expansion of $k(t,s)$ is known, then the method described above will produce the exact coefficients of the Taylor expansion of the solution $x(t)$. For the same reasons as those given for (3.9), the recursive formula (3.13) also produces the exact coefficients of a Taylor polynomial approximation to x(t) when the kernel k(t, s) and the function ξ(t) are sufficiently smooth but not necessarily analytic.

Under certain conditions, we can show the convergence of the power series with coefficients determined by (3.13). We define

$$L_0 := 0;$$
$$L_n := \sum_{i=0}^{n-1} \| K_{i,n-i-1} \|, \text{ for } n \geq 1$$

--- (3.14)

where the symbol $\|\cdot\|$ is used to indicate (depending on the context) a vector norm on $\mathrm{IR}^N$ and the induced matrix norm on $N \times N$ matrices.

We have the following:

<u>Theorem 3.1.</u> If each of the sets $\{\Xi_n : n = 0,1,2,...\}$, $\{L_n : n = 0,1,2,...\}$ is a set of coefficients of a convergent power series with strictly positive radius of convergence, then the series $\sum_{n=0}^{\infty} t^n X_n$, with $X_n$ given by (3.13), has a strictly positive radius of convergence.

<u>Proof:</u> We define the constants $C_n$, n=0, 1, 2, ... , by the recursive definition

$$C_0 := \|X_0\| \equiv \|\Xi_0\|;$$

$$C_n := \|\Xi_n\| + \sum_{l=0}^{n-1} L_{n-l} C_l \equiv \|\Xi_n\| + \sum_{l=0}^{n} L_{n-l} C_l \quad (\text{since } L_0 = 0), \text{ for } n \geq 1$$

--- (3.15)

It follows, by induction, that $\|X_n\| \leq C_n$ for all n=0, 1, 2, ... : if $\|X_l\| \leq C_l$ for $l = 0,1,2,...,n-1$, then we get, from (3.13),

$$\|X_n\| \leq \|\Xi\|_n + \sum_{l=0}^{n-1} \sum_{i=0}^{n-l-1} \frac{1}{n-i} \|K_{i,n-i-l-1}\| C_l \leq \|\Xi_n\| + \sum_{l=0}^{n-1} L_{n-l} C_l = C_n$$

--- (3.16)

Let $R_L$, $R_\Xi$ denote, respectively, the radii of convergence of power series with coefficients $L_n$, $\Xi_n$, n=0, 1, 2, ... .
In this proof, we shall use the symbol " $\underset{(FPS)}{=}$ " to denote equality in the sense of formal power series; for example, if we write

$$\left(\sum_{n=0}^{\infty} t^n a_n\right)\left(\sum_{n=0}^{\infty} t^n b_n\right) \underset{(FPS)}{=} \left(\sum_{n=0}^{\infty} t^n c_n\right),$$

then that equality will mean that $c_n = \sum_{l=0}^{n} a_{n-l} b_l$ for all n=0, 1, 2, ... , regardless of the convergence or divergence of the indicated power series. By contrast, a plain equality sign will denote equality of convergent power series with established radii of convergence.

Now, it follows from (3.15) that

$$\sum_{n=0}^{\infty} t^n \|C_n\| \stackrel{=}{_{(FPS)}} \left(\sum_{n=0}^{\infty} t^n \Xi_n\right) + \left(\sum_{n=0}^{\infty} t^n L_n\right)\left(\sum_{n=0}^{\infty} t^n C_n\right)$$

from which we obtain

$$\left(1 - \sum_{n=0}^{\infty} t^n L_n\right)\left(\sum_{n=0}^{\infty} t^n C_n\right) \stackrel{=}{_{(FPS)}} \left(\sum_{n=0}^{\infty} t^n \Xi_n\right)$$

--- (3.17)

The analytic function that is represented locally by the power series $1 - \sum_{n=0}^{\infty} t^n L_n$ takes a nonzero value (it becomes equal to 1) when t=0, therefore the function $\dfrac{1}{1 - \sum_{n=0}^{\infty} t^n L_n}$ can be expanded into a power series with strictly positive radius of convergence R',

$\dfrac{1}{1 - \sum_{n=0}^{\infty} t^n L_n} = \sum_{n=0}^{\infty} t^n M_n$ where the coefficients $M_n$ are determined by the algebraic

conditions $\left(\sum_{n=0}^{\infty} t^n M_n\right)\left(1 - \sum_{n=0}^{\infty} t^n L_n\right) \stackrel{=}{_{(FPS)}} 1$. Then (3.17) implies

$$\sum_{n=0}^{\infty} t^n C_n = \left(\sum_{n=0}^{\infty} t^n M_n\right)\left(\sum_{n=0}^{\infty} t^n \|\Xi_n\|\right)$$

--- (3.18)

where the product of the two power series on the right-hand side of (3.18) has radius of convergence at least $R_1 := \min(R_\Xi, R') > 0$, thus the coefficients $C_n$,

n=0, 1, 2, ... , and consequently also the coefficients $X_n$ , n=0, 1, 2, ... , are coefficients of convergent power series with strictly positive radii of convergence.  ///

The radius of convergence of the series $\sum_{n=0}^{\infty} t^n X_n$ , with $X_n$ given by (3.13), can be estimated as follows:

<u>Theorem 3.2.</u>  If the sets of coefficients $\{\Xi_n : n = 0,1,2,...\}$ , $\{L_n : n = 0,1,2,...\}$ satisfy

$|\Xi_n| \le \dfrac{M_0}{R^n}$ , $|L_n| \le \dfrac{M}{R^n}$ , for some $R > 0$, then the series $\sum_{n=0}^{\infty} t^n X_n$ has radius of convergence $R_X$ that satisfies $R_X \ge \dfrac{R}{1+RM}$ .

<u>Proof:</u>  We set

$$D_n := C_n R^n \qquad\qquad \text{--- (3.19)}$$

where the constants $C_n$ are given by (3.15). Then (3.15) and the inequalities $|\Xi_n| \le \dfrac{M_0}{R^n}$ and $|L_n| \le \dfrac{M}{R^n}$ imply

$$D_n \le M_0 + \sum_{l=0}^{n-1} MR D_l \qquad\qquad \text{--- (3.20)}$$

By a well known variant of Gronwall's inequality for discrete Volterra operators, (3.20) implies $D_n \le M_0 (1+MR)^n$ . For completeness, we include a very short proof: the statement is true for n=0 ; by induction, if $D_l \le M_0(1+MR)^l$ for $l = 0,1,2,...,n-1$, then

$$D_n \le M_0 + \sum_{l=0}^{n-1} MR M_0 (1+MR)^l = M_0 \left[ 1 + MR \dfrac{(1+MR)^n - 1}{MR} \right] = M_0 (1+MR)^n$$

$$\text{--- (3.21)}$$

We look at the coefficients $E_n := M_0 \frac{(1+MR)^n}{R^n}$, and we note that $0 \leq \|X_n\| \leq C_n \leq E_n$. A ratio test shows that a power series with coefficients $E_n$ would have radius of convergence $R_E = \frac{R}{1+RM}$, and thus, by a comparison argument, the power series with coefficients $X_n$ has radius of convergence not less than $R_E$. ///

## 4. Systems with Abel-type singularities.

We consider the N-dimensional system

$$x(t) = \xi(t) + \int_0^t f(t, s, x(s))ds$$

--- (4.1)

where the function f has the following form:

$$f(t,s,x) = a(t,s,x) + (t-s)^{-\alpha} b(t,s,x), \quad t > s, \quad 0 < \alpha < 1;$$

$$a(t,s,x) = \sum_{i=0}^{\infty} \sum_{j=0}^{\infty} \sum_{k \geq 0} t^i s^j x^k A_{ijk};$$

$$b(t,s,x) = \sum_{i=0}^{\infty} \sum_{j=0}^{\infty} \sum_{k \geq 0} t^i s^j x^k B_{ijk}$$

--- (4.2)

It turns out that it is easier to start with the case in which the exponent $\alpha$ is irrational. In that case, we take

$$\xi(t) = \sum_{r=0}^{\infty} \sum_{i > r\alpha - 1} t^{i - r\alpha} \Xi_{ri}$$

--- (4.3)

and we seek a solution of (4.1) in the form

$$x(t) = \sum_{r=0}^{\infty} \sum_{i > r\alpha - 1} t^{i - r\alpha} X_{ri}$$

--- (4.4)

The reason for taking x(t) in the form of a series in the functions $\{t^{i-r\alpha} : r \geq 0, i > r\alpha - 1\}$ is that the set of these functions is the smallest set of functions that remains closed under the operations of substituting a series of the form of (4.4) into the series expansions (4.2) for a(t, s, x(s)) and b(t, s, x(s)) and then evaluating the integral
$\int_0^t \{a(t,s,x(s)) + (t-s)^{-\alpha} b(t,s,x(s))\} ds$. The closure properties of the set
$\{t^{i-r\alpha} : r \geq 0, i > r\alpha - 1\}$ will become evident after we present the relevant calculations below. The calculations below are predicated on the validity of a series solution of the indicated type with positive radius of convergence.
In addition to the notation and conventions of section 3, we shall use the following:

$$Y_{\rho l} := \{X_{ri} : r \leq \rho, i \leq l\};$$

$$\mathbf{r} := (r_1, r_2, ..., r_N); |\mathbf{k}| := k_1 + k_2 + ... + k_N;$$

$$X_{\mathbf{r}i} := \{X_{r_j i_j(j)} : j = 1, 2, ..., N\}; (X_{\mathbf{r}i})^{\mathbf{k}} := \prod_{j=1}^{N} \left(X_{r_j i_j(j)}\right)^{k_j}$$

--- (4.5)

A power $(x(t))^{\mathbf{k}}$, when $x(t)$ is given by (4.4), is evaluated as follows:

$$(x(t))^{\mathbf{k}} = \prod_{j=1}^{N} \left( \sum_{r=0}^{\infty} \sum_{i > r\alpha - 1} t^{i - r\alpha} X_{ri(j)} \right)^{k_j} =$$

$$= \sum_{\mathbf{r} \geq 0} \sum_{\mathbf{i} : i_j > r_j \alpha - 1, \forall j} t^{\mathbf{i} \bullet \mathbf{k} - \alpha \mathbf{r} \bullet \mathbf{k}} \prod_{j=1}^{N} \left(X_{r_j i_j(j)}\right)^{k_j} =$$

$$= \sum_{\rho=0}^{\infty} \sum_{l > \rho\alpha - |\mathbf{k}|} \sum_{\mathbf{r} : \mathbf{r} \bullet \mathbf{k} = \rho} \sum_{\substack{\mathbf{i} : i_j > r_j \alpha - 1, \forall j; \\ \mathbf{i} \bullet \mathbf{k} = l}} t^{l - \rho\alpha} \prod_{j=1}^{N} \left(X_{r_j i_j(j)}\right)^{k_j} =$$

$$= \sum_{\rho=0}^{\infty} \sum_{l > \rho\alpha - |\mathbf{k}|} \sum_{\mathbf{r} : \mathbf{r} \bullet \mathbf{k} = \rho} \sum_{\substack{\mathbf{i} : i_j > r_j \alpha - 1, \forall j; \\ \mathbf{i} \bullet \mathbf{k} = l}} t^{l - \rho\alpha} (X_{\mathbf{r}i})^{\mathbf{k}}$$

--- (4.6)

We use the notation

$$Z_{\mathbf{k}}(Y_{\rho l}) := \sum_{\mathbf{r}:\, \mathbf{r} \bullet \mathbf{k} = \rho} \sum_{\substack{\mathbf{i}:\, i_j > r_j \alpha - 1,\, \forall j; \\ \mathbf{i} \bullet \mathbf{k} = l}} (X_{\mathbf{ri}})^{\mathbf{k}}$$

--- (4.7)

Then we have

$$(x(t))^{\mathbf{k}} = \sum_{\rho=0}^{\infty} \sum_{l > \rho\alpha - |\mathbf{k}|} t^{l-\rho\alpha} Z_{\mathbf{k}}(Y_{\rho l})$$

--- (4.8)

For later reference, we write down these two integrals:

$$\int_0^t t^i s^j s^{l-r\alpha} ds = \frac{t^{i+j+l-r\alpha+1}}{j+l-r\alpha+1};$$

$$\int_0^t (t-s)^{-\alpha} t^i s^j s^{l-r\alpha} ds = B(1-\alpha, j+l-r\alpha+1) t^{i+j+l-(r+1)\alpha+1}$$

--- (4.9)

where $B(\mu,\nu)$ is the standard beta function, i.e. $B(\mu,\nu) = \dfrac{\Gamma(\mu)\Gamma(\nu)}{\Gamma(\mu+\nu)}$ for $\mu > 0$, $\nu > 0$.
(The beta function $B(\mu,\nu)$ is, of course, unrelated to the Taylor coefficients $B_{ij\mathbf{k}}$.) We shall denote by $\chi$ the truth function, i.e., if $P$ is a logical statement, then $\chi(P) = 1$ if $P$ is true, and $\chi(P) = 0$ if $P$ is false.

After substituting the series expansions into (4.1), and using also the integral formulae (4.9), we obtain the following recursive equation for the coefficients $X_{\rho n}$:

$$X_{\rho n} = \Xi_{\rho n} + \sum_{\rho\alpha-1 < l \leq n-1} \sum_{i=0}^{n-l-1} \sum_{\mathbf{k} \in I_l \cap I_\rho} \frac{1}{n-i-\rho\alpha} Z_{\mathbf{k}}(Y_{\rho l}) A_{i, n-i-l-1, \mathbf{k}} +$$

$$+ \chi(\rho \geq 1) \sum_{(\rho-1)\alpha-1 < l \leq n-1} \sum_{i=0}^{n-l-1} \sum_{\mathbf{k} \in I_l \cap I_{\rho-1}} B(1-\alpha, n-i-(\rho-1)\alpha) Z_{\mathbf{k}}(Y_{\rho-1, l}) B_{i, n-i-l-1, \mathbf{k}}$$

--- (4.10)

Now we consider the case in which $\alpha$ is rational, say $\alpha = \dfrac{p}{q}$, with $0 < p < q$ and $p, q$ relatively prime.

We assume that the function $\xi(t)$ has the expansion

$$\xi(t) = \sum_{r=0}^{q-1} \sum_{j > r\alpha - 1} t^{j - r\alpha} \Xi_{rj}$$

--- (4.11)

In this case we seek a solution $x(t)$ in the form

$$x(t) = \sum_{r=0}^{q-1} \sum_{i > r\alpha - 1} t^{i - r\alpha} X_{ri}$$

--- (4.12)

We note that the exponent $i + j + l - (r+1)\alpha + 1$, that appears on the right-hand side of each equation in (4.9), is non-integer for $0 \le r \le q - 2$, but it becomes equal to $i + j + l - p + 1$ when $r = q - 1$.

Then, substitution of the series expansions into the integral equation (4.1) gives

$$X_{\rho n} = \Xi_{\rho n} + \sum_{\rho\alpha - 1 < l \le n - 1} \sum_{i=0}^{n-l-1} \sum_{\mathbf{k} \in I_l \cap I_\rho} \frac{1}{n - i - \rho\alpha} Z_{\mathbf{k}}(Y_{\rho l}) A_{i, n-i-l-1, \mathbf{k}} +$$

$$+ \chi(\rho \ge 1) \sum_{(\rho-1)\alpha - 1 < l \le n - 1} \sum_{i=0}^{n-l-1} \sum_{\mathbf{k} \in I_l \cap I_{\rho-1}} B(1 - \alpha, n - i - (\rho - 1)\alpha) \cdot$$

$$\cdot Z_{\mathbf{k}}(Y_{\rho-1, l}) B_{i, n-i-l-1, \mathbf{k}} +$$

$$+ \chi(\rho = q) \sum_{(\rho-1)\alpha - 1 < l \le n+p-1} \sum_{i=0}^{n+p-l-1} B(1 - \alpha, n - i - (\rho - 1)\alpha) Z_{\mathbf{k}}(Y_{\rho-1, l}) B_{i, n+p-i-l-1, \mathbf{k}}$$

--- (4.13)

Another method for obtaining the coefficients of $x(t)$ when $\alpha$ is rational would be to first use a formal expression

$$x(t) = \sum_{r=0}^{\infty} \sum_{i > r\alpha - 1} t^{i - r\alpha} W_{ri}$$

--- (4.14)

and calculate the coefficients $W_{ri}$ by the same method as for the case of irrational $\alpha$:

$$W_{\rho n} = \Xi_{\rho n} + \sum_{\rho\alpha - 1 < l \leq n - 1} \sum_{i=0}^{n-l-1} \sum_{\mathbf{k} \in I_l \cap I_\rho} \frac{1}{n - i - \rho\alpha} Z_\mathbf{k}(Y^W_{\rho l}) A_{i, n-i-l-1, \mathbf{k}} +$$

$$+ \chi(\rho \geq 1) \sum_{(\rho-1)\alpha - 1 < l \leq n - 1} \sum_{i=0}^{n-l-1} \sum_{\mathbf{k} \in I_l \cap I_{\rho-1}} B(1 - \alpha, n - i - (\rho - 1)\alpha) \cdot$$

$$\cdot Z_\mathbf{k}(Y^W_{\rho-1, l}) B_{i, n-i-l-1, \mathbf{k}}$$

--- (4.15)

where $Y_{\rho l}^W := \{W_{ri} : r \leq \rho, i \leq l\}$, $Z_\mathbf{k}(Y_{\rho l}^W) := \sum_{\mathbf{r} : \mathbf{r} \bullet \mathbf{k} = \rho} \sum_{\substack{\mathbf{i} : i_j > r_j \alpha - 1, \forall j; \\ \mathbf{i} \bullet \mathbf{k} = l}} (W_{\mathbf{ri}})^\mathbf{k}$.

We observe that, when $\alpha = \dfrac{p}{q}$, with $0 < p < q$ and $p, q$ relatively prime, each nonnegative integer r is uniquely representable as $r = mq + \rho$ with $\rho \in \{0, 1, ..., q - 1\}$, and the corresponding exponent of t, in the expansion $\sum_{r=0}^{\infty} \sum_{i > r\alpha - 1} t^{i - r\alpha} W_{ri}$ becomes $i - r\alpha = i - mq\alpha - \rho\alpha = i - mp - \rho\alpha$, and consequently, if we want that exponent to be equal to $n - \rho\alpha$, we must have i=n+mp. Therefore, the coefficients $X_{ri}$ in (4.12) are related to the coefficients $W_{ri}$ via

$$X_{\rho n} = \sum_{m=0}^{\infty} W_{\rho + mq, n + mp}, \quad \rho = 0, 1, 2, ..., q - 1$$

--- (4.16)

As in the previous sections, also for the problems of this section, the recursive relation can be specialized to the case of linear Volterra systems. We present this case below.

We consider the system of Volterra integral equations

$$x(t) = \xi(t) + \int_0^t k(t,s)x(s)ds$$

--- (4.17)

where the dimensionality of $x(t)$ and $k(t,s)$ is as in section 2 above. Here, however, the kernel $k(t,s)$ has the expression

$$k(t,s) = a(t,s) + (t-s)^{-\alpha} b(t,s), \, t > s, \, 0 < \alpha < 1$$

--- (4.18)

We assume that $a(t,s)$ and $b(t,s)$ have the expansions

$$a(t,s) = \sum_{i=0}^{\infty} \sum_{j=0}^{\infty} t^i s^j A_{ij}, \, b(t,s) = \sum_{i=0}^{\infty} \sum_{j=0}^{\infty} t^i s^j B_{ij}$$

--- (4.19)

As before, we start with the case in which $\alpha$ is irrational.

We assume that $\xi(t)$ can be represented as

$$\xi(t) = \sum_{r=0}^{\infty} \sum_{i > r\alpha - 1} t^{i-r\alpha} \Xi_{ri}$$

--- (4.20)

and we seek a solution in the form

$$x(t) = \sum_{r=0}^{\infty} \sum_{i>r\alpha-1} t^{i-r\alpha} X_{ri}$$

--- (4.21)

By using the integral formulae (4.9), and substituting the series expressions into the integral equation (4.17), we find

$$X_{rn} = \Xi_{rn} + \sum_{r\alpha-1<l\leq n-1} \sum_{i=0}^{n-l-1} \frac{1}{n-i-r\alpha} A_{i,n-i-l-1} X_{rl} +$$

$$+ \chi(r \geq 1) \sum_{(r-1)\alpha-1<l\leq n-1} \sum_{i=0}^{n-l-1} B(1-\alpha, n-i-(r-1)\alpha) B_{i,n-i-l-1} X_{r-1,l}$$

--- (4.22)

Eq. (4.22) is a simple recursive formula for the unknown coefficients of the solution $x(t)$. It is also a discrete version of an integro-differential Volterra equation in two variables, i.e. (4.22) is a finite-difference equation in the variable r and a discrete Volterra equation in the variable n.

Now we look at the case in which $\alpha$ is rational, say $\alpha = \frac{p}{q}$, with $0 < p < q$ and $p,q$ relatively prime.

We assume that the function $\xi(t)$ has the expansion

$$\xi(t) = \sum_{r=0}^{q-1} \sum_{j>r\alpha-1} t^{j-r\alpha} \Xi_{rj}$$

--- (4.23)

In this case we take a solution $x(t)$ in the form

$$x(t) = \sum_{r=0}^{q-1} \sum_{i>r\alpha-1} t^{i-r\alpha} X_{ri}$$

--- (4.24)

Then, substitution of the series expansions into the integral equation (4.17) gives

$$X_{rn} = \Xi_{rn} + \sum_{r\alpha-1<l\leq n-1} \sum_{i=0}^{n-l-1} \frac{1}{n-i-r\alpha} A_{i,n-i-l-1} X_{rl} +$$

$$+ \chi(1 \leq r \leq q-1) \sum_{(r-1)\alpha-1<l\leq n-1} \sum_{i=0}^{n-l-1} B(1-\alpha, n-i-(r-1)\alpha) B_{i,n-i-l-1} X_{r-1,l} +$$

$$+ \chi(r=q) \sum_{(r-1)\alpha-1<l\leq n+p-1} \sum_{i=0}^{n+p-l-1} B(1-\alpha, n-i-(r-1)\alpha+1) B_{i,n+p-i-l-1} X_{r-1,l}$$

--- (4.25)

Eq. (4.25) is still a recursive formula for the coefficients $X_{rn}$ but it has a peculiar feature: for $r = q$, the calculation of $X_{rn}$ requires knowledge of $X_{r-1,l}$ for values of $l$ greater than n, and therefore it is not a discrete Volterra equation in the variable n.

As in the general case of the nonlinear Volterra system (4.1), there is also a second method to get the coefficients of x(t) in the linear case.
First, we denote by $Z_{rn}$ the formal coefficients obtained by using (4.22) but with rational $\alpha$:

$$Z_{rn} = \Xi_{rn} + \sum_{r\alpha-1<l\leq n-1} \sum_{i=0}^{n-l-1} \frac{1}{n-i-r\alpha} A_{i,n-i-l-1} Z_{rl} +$$

$$+ \chi(r \geq 1) \sum_{(r-1)\alpha-1<l\leq n-1} \sum_{i=0}^{n-l-1} B(1-\alpha, n-i-(r-1)\alpha) B_{i,n-i-l-1} Z_{r-1,l}$$

--- (4.26)

Since $p = \alpha q$ with p,q relatively prime, the terms in the expansion for x(t) that contain $t^{n-r\alpha}$, $r = 0,1,2,...,q-1$, are precisely all the terms with coefficients $Z_{r+kp,n+kp}$, $k = 0,1,2,...$, thus

$$X_{rn} = \sum_{k=0}^{\infty} Z_{r+kq, n+kp}$$

--- (4.27)

In this second method, we obtain a discrete Volterra system for the calculation of the coefficients $Z_{rn}$, but then the evaluation of $X_{rn}$ requires an infinite number of coefficients $Z_{r'n'}$.

Example 4.1.  The solution to the simple scalar linear Abel integral equation of the second kind,

$$x(t) = \xi(t) + c\int_0^t (t-s)^{-\alpha} x(s)\,ds, \; 0 < \alpha < 1$$

--- (4.28)

has been reported in [PM], and it has been obtained by methods different from those we have used in the present paper. Also, for rational values of $\alpha$, the same Abel equation of the second kind has been solved in [BNR]. The solution in [PM] holds for both, rational and irrational $\alpha$, and for this reason we use the form of the solution given in [PM]. This equation is, of course, a particular case of (4.17) with $a(t,s) \equiv 0$, $b(t,s) \equiv c = $ constant (so that, in particular, $A_{ij} = 0$ for all i and j, $B_{00} = c$, and $B_{ij} = 0$ if $(i,j) \neq (0,0)$ ).

After we make some simple changes to fit the notation of the present paper, the solution in [PM] amounts to

$$x(t) = \xi(t) + \int_0^t R(t-s)x(s)\,ds;$$

$$R(s) = \sum_{n=0}^{\infty} P(n,c,\alpha) t^{n-(n+1)\alpha}; \; P(n,c,\alpha) := \frac{c^{n+1}(\Gamma(1-\alpha))^{n+1}}{\Gamma((n+1)(1-\alpha))}$$

--- (4.29)

For the purposes of this example, we assume that the function $\xi$ is analytic, and has the expansion

$$\xi(t) = \sum_{i=0}^{\infty} t^i \Xi_{0i}$$

--- (4.30)

thus $\Xi_{ri} = 0$ for $r \geq 1$.
Now, we can calculate

$$\int_0^t R(t-s)\xi(s)\,ds = \sum_{r=1}^{\infty} \sum_{n=r}^{\infty} t^{n-r\alpha} B(r(1-\alpha), n-r-1) P(r-1, \alpha, c) \Xi_{0,n-r}$$

--- (4.31)

and therefore, in our terminology,

$$X_{0n} = \Xi_{0n} ;$$
$$X_{rn} = \chi(n \geq r) B(r(1-\alpha), n-r-1) P(r-1, \alpha, c) \Xi_{0,n-r}, \text{ for } r \geq 1$$

--- (4.32)

Thus, in this example, x(t) does indeed have an expansion of the type (4.21), and consequently those coefficients are the same with the solution of the recursive relationship (4.22).
It can also be directly verified, by using (4.22), that the expression (4.32) is the unique solution of (4.22) for the case of this example. We note that, in this case, (4.22) becomes

$$X_{0n} = \Xi_{0n} ;$$
$$X_{rn} = cB(1-\alpha, n-(r-1)\alpha) X_{r-1,n-1}$$

--- (4.33)

and the solution of (4.33) is

$$X_{rn} = c^r \Xi_{0,n-r} \prod_{j=1}^{r} B(1-\alpha, n-(j-1)-(r-j)\alpha) =$$
$$= c^r \Xi_{0,n-r} \prod_{j=1}^{r} \frac{\Gamma(1-\alpha)\Gamma(n-(j-1)-(r-j)\alpha)}{\Gamma(n-(j-2)-(r-(j-1))\alpha)} = c^r \Xi_{0,n-r} [\Gamma(1-\alpha)]^r \frac{\Gamma(1+n-r)}{\Gamma(1+n-r\alpha)} ,$$
for $r \geq 1$

--- (4.34)

where, in the last step in (4.34), we have used suitable cancellations.
On the other hand, the solution given by (4.32) becomes, for $r \geq 1$,

$$X_{rn} = c^r \Xi_{0,n-r} \frac{\Gamma(r(1-\alpha))\Gamma(1+n-r)}{\Gamma(1+n-r\alpha)} \frac{[\Gamma(1-\alpha)]^r}{\Gamma(r(1-\alpha))} =$$

$$= c^r \Xi_{0,n-r} [\Gamma(1-\alpha)]^r \frac{\Gamma(1+n-r)}{\Gamma(1+n-r\alpha)}$$

--- (4.35)

thus we have directly verified that the solution obtained by using (4.22) coincides with the solution in [PM]. ///

Next, we analyze the convergence of the series with coefficients given by (4.22), in the case in which $\alpha$ is rational. We take an arbitrary, but fixed throughout this analysis, number $\delta > 0$, and we examine the series

$$S(t) :=_{(FPS)} \sum_{r=0}^{\infty} \sum_{i:\, i > r\alpha - 1} t^{i-r\alpha} X_{ri} \text{ for } t \geq \delta$$

--- (4.36)

By formally rearranging the order of summation in (4.36), we obtain

$$S(t) =_{(FPS)} \sum_{i=0}^{\infty} \sum_{r:\, r\alpha < i+1} t^{i-r\alpha} X_{ri}$$

--- (4.37)

For $t \geq \delta$, we have

$$\sum_{r:\, r\alpha < i+1} t^{i-\rho\alpha} X_{ri} \leq t^i \sum_{r:\, r\alpha < i+1} \delta^{-r\alpha} X_{ri}$$

--- (4.38)

We set

$$\mathbf{X}_i := \sum_{r:\, r\alpha < i+1} \delta^{-r\alpha} \|X_{ri}\| \,;\, \mathbf{\Xi}_i := \sum_{r:\, r\alpha < i+1} \delta^{-r\alpha} \|\Xi_{ri}\|$$

--- (4.39)

We have

<u>Lemma 4.1.</u> If $\alpha = \dfrac{p}{q}$ with p, q relatively prime positive integers and p<q, then, for $0 \le i \le n - l - 1$ and $r\alpha - 1 < l \le n - 1$, we have
$$\dfrac{1}{n - i - r\alpha} \le q \,,\; B(1 - \alpha, n - i - r\alpha) \le B(1 - \alpha, \tfrac{1}{q}).$$

<u>Proof:</u> There are two possibilities to consider: $r\alpha$ is an integer, or it is not an integer. If $r\alpha$ is an integer, then the condition $r\alpha - 1 < l$ amounts to $l \ge r\alpha$; the minimum value of $n - i - r\alpha$, for fixed $l$, is $l + 1 - r\alpha$, and the minimum value of $l + 1 - r\alpha$ is 1, thus $\dfrac{1}{n - i - r\alpha} \le 1$. For the inequality with the beta function, we use the fact that $B(1-\alpha, \beta)$ is a decreasing function of $\beta$, and we conclude that $B(1-\alpha, n - i - r\alpha) \le B(1-\alpha, 1)$. When $r\alpha$ is not an integer, a similar argument as above shows that the minimum value of $n - i - r\alpha$, under the stated restrictions on i and n, is $\text{fract}(r\alpha)$, the fractional part of $r\alpha$, defined as $\text{fract}(r\alpha) := r\alpha - \text{floor}(r\alpha)$, where the floor function (a.k.a. greatest integer function) of a real number u is defined in the standard way, i.e. it is the greatest among all integers that are $\le u$. The minimum possible value of $\text{fract}(r\alpha)$ is $\dfrac{1}{q}$, and consequently $\dfrac{1}{n - i - r\alpha} \le q$ and $B(1-\alpha, n - i - r\alpha) \le B(1-\alpha, \tfrac{1}{q})$. ///

We define

$$L_m := \sum_{i=0}^{m-1} \|A_{i, m-i-1}\| \,,\; M_m := \sum_{i=0}^{m-1} \|B_{i, m-i-1}\| \,,\; \text{for } m \ge 1;$$
$$L_0 := 0 \,;\, M_0 := 0$$

--- (4.40)

Then we obtain from (4.22)

$$\|X_{rn}\| \leq \|\Xi_{rn}\| + \sum_{l:\, r\alpha-1<l\leq n} qL_{n-l}\|X_{rl}\| +$$

$$+ \chi(r \geq 1) \sum_{l:\, (r-1)\alpha-1<l\leq n} B(1-\alpha, \frac{1}{q})M_{n-l}\|X_{r-1,l}\|$$

--- (4.41)

from which it follows that

$$\sum_{r:\, r\alpha-1<n} \delta^{-r\alpha}\|X_{rn}\| \leq \sum_{r:\, r\alpha-1<n} \delta^{-r\alpha}\|\Xi_{rn}\| +$$

$$+ \sum_{r:\, r\alpha-1<n} \sum_{l:\, r\alpha-1<l\leq n} \delta^{-r\alpha} qL_{n-l}\|X_{rl}\| +$$

$$+ \sum_{r:\, r\alpha-1<n} \chi(r \geq 1) \sum_{l:\, (r-1)\alpha-1<l\leq n} \delta^{-r\alpha}B(1-\alpha, \frac{1}{q})M_{n-l}\|X_{r-1,l}\|$$

--- (4.42)

We note that $\{r: r\alpha - 1 < n\} \supseteq \{r: (r-1)\alpha - 1 < n\}$, and therefore (4.42) implies

$$\sum_{r:\, r\alpha-1<n} \delta^{-r\alpha}\|X_{rn}\| \leq \sum_{r:\, r\alpha-1<n} \delta^{-r\alpha}\|\Xi_{rn}\| +$$

$$+ \sum_{r:\, r\alpha-1<n} \sum_{l:\, r\alpha-1<l\leq n} \delta^{-r\alpha} qL_{n-l}\|X_{rl}\| +$$

$$+ \sum_{r:\, (r-1)\alpha-1<n} \chi(r \geq 1) \sum_{l:\, (r-1)\alpha-1<l\leq n} \delta^{-r\alpha}B(1-\alpha, \frac{1}{q})M_{n-l}\|X_{r-1,l}\|$$

--- (4.43)

By changing the order of summations in (4.43), we obtain

$$\sum_{r:\, r\alpha-1<n} \delta^{-r\alpha} \|X_{rn}\| \leq \sum_{r:\, r\alpha-1<n} \delta^{-r\alpha} \|\Xi_{rn}\| +$$

$$+ \sum_{l=0}^{n} \sum_{r:\, r\alpha-1<l} \delta^{-r\alpha} q L_{n-l} \|X_{rl}\| +$$

$$+ \sum_{l=0}^{n} \sum_{r:\, r\alpha-1<l} \chi(r \geq 1)\, \delta^{-\alpha} \delta^{-(r-1)\alpha} B\!\left(1-\alpha, \frac{1}{q}\right) M_{n-l} \|X_{r-1,l}\|$$

--- (4.44)

and, in view of the definitions in (4.39), we obtain

$$\mathbf{X}_n \leq \mathbf{\Xi}_n + \sum_{l=0}^{n} \left[ q L_{n-l} + \delta^{-\alpha} B\!\left(1-\alpha, \frac{1}{q}\right) M_{n-l} \right] \mathbf{X}_l =$$

$$= \mathbf{\Xi}_n + \sum_{l=0}^{n-1} \left[ q L_{n-l} + \delta^{-\alpha} B\!\left(1-\alpha, \frac{1}{q}\right) M_{n-l} \right] \mathbf{X}_l$$

--- (4.45)

The inequality (4.45) is of the same type as (3.16), and it is clear that the same methods of proof as in theorems (3.1) and (3.2) also apply to this case. Because of the similarity in the methods of proof, we state the following two theorems without proof.

Theorem 4.1. We assume that $\alpha$ is rational of the form described in lemma 4.1. If each of the sets $\{\mathbf{\Xi}_n : n = 0,1,2,...\}$, $\{L_n : n = 0,1,2,...\}$, $\{M_n : n = 0,1,2,...\}$ is a set of coefficients of a convergent power series with strictly positive radius of convergence, then the series $\sum_{n=0}^{\infty} t^n \mathbf{X}_n$, with $\mathbf{X}_n$ given by (4.39), has a strictly positive radius of convergence. ///

Theorem 4.2. We assume that $\alpha$ is rational of the form described in lemma 4.1. If the sets of coefficients $\{\mathbf{\Xi}_n : n = 0,1,2,...\}$, $\{L_n : n = 0,1,2,...\}$, $\{M_n : n = 0,1,2,...\}$ satisfy

$$|\Xi_n| \leq \frac{C_0}{R^n}, \ |L_n| \leq \frac{C_L}{R^n}, \ |M_n| \leq \frac{C_M}{R^n}, \text{ for some } R > 0, \text{ then the series}$$

$$S_1(t) \underset{(FPS)}{:=} \sum_{n=0}^{\infty} t^n \mathbf{X}_n \text{ has radius of convergence } R_X \text{ that satisfies } R_X \geq \frac{R}{1 + RC_\delta} \text{ where}$$

$$C_\delta := qC_L + \delta^{-\alpha} B\left(1 - \alpha, \frac{1}{q}\right) C_M. \ ///$$

When the series $S_1(t)$ (defined in the statement of Theorem 4.2 above) has radius of convergence $R_X$, then the series $S(t)$ (defined in (4.36)) converges absolutely for $\delta \leq t < R_X$; this follows from

$$\sum_{i=0}^{\infty} \sum_{r: \, r\alpha - 1 < i} t^{i - r\alpha} |X_{ri}| \leq \sum_{i=0}^{\infty} \sum_{r: \, r\alpha - 1 < i} t^i \delta^{-r\alpha} |X_{ri}| = \sum_{i=0}^{\infty} t^i \mathbf{X}_i$$

--- (4.46)

Consequently, the series $S(t)$ represents a solution of the integral equation (4.17) for $\delta \leq t < R_X$, for every $\delta > 0$. Of course, the value of $R_X$, estimated in theorem 4.2, depends on $\delta$.

## 5. Systems with logarithmic singularities.

We consider the nonlinear N-dimensional Volterra system

$$x(t) = \xi(t) + \int_0^t f(t,s,x(s))ds$$

--- (5.1)

with

$$f(t,s,x) = a(t,s,x) + \ln(t-s)b(t,s,x), \text{ for } t > s;$$

$$a(t,s,x) = \sum_{i=0}^{\infty} \sum_{j=0}^{\infty} \sum_{k \geq 0} t^i s^j x^k A_{ijk};$$

$$b(t,s,x) = \sum_{i=0}^{\infty} \sum_{j=0}^{\infty} \sum_{k \geq 0} t^i s^j x^k B_{ijk}$$

--- (5.2)

We take

$$\xi(t) = \sum_{r=0}^{\infty} \sum_{i=0}^{\infty} (\ln t)^r t^i \Xi_{ri}$$

--- (5.3)

and we seek a solution $x(t)$ of the form

$$x(t) = \sum_{r=0}^{\infty} \sum_{i=0}^{\infty} (\ln t)^r t^i X_{ri}$$

--- (5.4)

As in the previous section, we are concerned with the evaluation of the coefficients $X_{ri}$, assuming the validity of the indicated series solution with a strictly positive radius of convergence. The system of functions $\{(\ln t)^r t^i : r \geq 0, i \geq 0\}$ is the smallest set of functions that remains invariant after substitution of a series of the type of (5.4) into the series expansions for a(t, s, x(s)) and b(t, s, x(s)), and then evaluation of the integral $\int_0^t \{a(t,s,x(s)) + (\ln(t-s))b(t,s,x(s))\}ds$, i.e. if $x(t)$ is represented as a series in the functions $\{(\ln t)^r t^i : r \geq 0, i \geq 0\}$, then $\int_0^t \{a(t,s,x(s)) + (\ln(t-s))b(t,s,x(s))\}ds$ is also representable as a series in the same functions.

We shall need the following integrals:

$$L_{qr} := \int_0^1 \sigma^q (\ln \sigma)^r d\sigma, \quad M_{qr} := \int_0^1 \sigma^q (\ln \sigma)^r \ln(1-\sigma) d\sigma$$

--- (5.5)

These two expressions, $L_{qr}$ and $M_{qr}$, are useful for evaluating the integrals that arise from substitution of the relevant series into the integral equation. Indeed, by using the change of variables $\sigma := \frac{s}{t}$, we find

$$\int_0^t t^i s^j s^l (\ln s)^r ds = t^{i+j+l+1} \int_0^1 \sigma^{j+l} (\ln t + \ln \sigma)^r d\sigma =$$

$$= t^{i+j+l+1} \int_0^1 \sigma^{j+l} \sum_{\mu=0}^r \binom{r}{\mu} (\ln t)^\mu (\ln \sigma)^{r-\mu} d\sigma =$$

$$= t^{i+j+l+1} \sum_{\mu=0}^r \binom{r}{\mu} (\ln t)^\mu L_{j+l, r-\mu};$$

$$\int_0^t t^i s^j s^l (\ln s)^r \ln(t-s) ds = t^{i+j+l+1} \int_0^1 \sigma^{j+l} (\ln t + \ln \sigma)^r (\ln t)\ln(1-\sigma) d\sigma =$$

$$= t^{i+j+l+1} \int_0^1 \sigma^{j+l} \sum_{\mu=0}^r \binom{r}{\mu} (\ln t)^{\mu+1} \ln(1-\sigma) (\ln \sigma)^{r-\mu} d\sigma =$$

$$= t^{i+j+l+1} \sum_{\mu=0}^r \binom{r}{\mu} (\ln t)^{\mu+1} M_{j+l, r-\mu}$$

--- (5.6)

The two integrals in (5.5) are evaluated as follows:

Using elementary integration techniques, with the substitutions $u = -\ln \sigma$, $v = (q+1)u$, we find

$$L_{qr} \equiv \int_0^1 \sigma^q (\ln \sigma)^r \, d\sigma = (-1)^r \int_0^\infty u^r \exp(-(q+1)u) \, du =$$

$$= \frac{(-1)^r}{(q+1)^{r+1}} \int_0^\infty v^r \exp(-v) \, dv = \frac{(-1)^r}{(q+1)^{r+1}} \Gamma(r+1) = \frac{(-1)^r (r!)}{(q+1)^{r+1}}$$

--- (5.7)

For the quantities $M_{qr}$ we use a series expansion for $\ln(1-\sigma)$, namely

$$\ln(1-\sigma) = \sum_{\lambda=0}^{\infty} \frac{\sigma^\lambda}{\lambda+1}, \text{ valid for } -1 \leq \sigma < 1, \text{ and we calculate, formally at first,}$$

$$M_{qr} \equiv \int_0^1 \sigma^q (\ln \sigma)^r \ln(1-\sigma) \, d\sigma = \sum_{\lambda=0}^{\infty} \int_0^1 \frac{\sigma^{q+\lambda+1}}{\lambda+1} (\ln \sigma)^r \, d\sigma =$$

$$= \sum_{\lambda=0}^{\infty} \frac{1}{\lambda+1} L_{q+\lambda+1,r} = \sum_{\lambda=0}^{\infty} \frac{(-1)^r (r!)}{(\lambda+1)(q+\lambda+1)^{r+1}}$$

--- (5.8)

and the validity of the above formal calculation now follows from the convergence of the last series in (5.8) and Abel's lemma.

Eq. (5.8) is the best we can do about $M_{qr}$ for general values of q and r; in the particular case $r = 0$, it is possible to express the term $M_{q0}$ in closed form. By setting $u = -\ln(1-\sigma)$, we have

$$M_{q0} \equiv \int_0^1 \sigma^q \ln(1-\sigma) \, d\sigma = \int_0^\infty u \left(1 - e^{-u}\right)^q e^{-u} \, du = \int_0^\infty \sum_{\mu=0}^{q} (-1)^\mu \binom{q}{\mu} u \, e^{-(\mu+1)u} \, du =$$

$$= \sum_{\mu=0}^{q} \frac{(-1)^\mu}{(\mu+1)^2} \binom{q}{\mu}$$

--- (5.9)

Now, we return to the Volterra system (5.1).
By substituting the series expansions into (5.1), we obtain the equations

$$\sum_{\rho=0}^{\infty}\sum_{\lambda=0}^{\infty}(\ln t)^{\rho}t^{\lambda}X_{\rho\lambda} = \sum_{\rho=0}^{\infty}\sum_{\lambda=0}^{\infty}(\ln t)^{\rho}t^{\lambda}\Xi_{\rho\lambda} +$$

$$+\sum_{i\geq 0,\ j\geq 0,\ \mathbf{k}\geq 0}\sum_{\rho\geq\mu}\sum_{\lambda=0}^{\infty}\sum_{\mathbf{r}:\ \mathbf{r}\bullet\mathbf{k}=\rho}\sum_{\mathbf{l}:\ \mathbf{l}\bullet\mathbf{k}=\lambda} t^{i+j+\lambda+1}(\ln t)^{\mu}\binom{\rho}{\mu}L_{j+1,\rho-\mu}(X_{\mathbf{r}\mathbf{l}})^{\mathbf{k}}A_{ij\mathbf{k}} +$$

$$+\sum_{i\geq 0,\ j\geq 0,\ \mathbf{k}\geq 0}\sum_{\rho\geq\mu}\sum_{\lambda=0}^{\infty}\sum_{\mathbf{r}:\ \mathbf{r}\bullet\mathbf{k}=\rho}\sum_{\mathbf{l}:\ \mathbf{l}\bullet\mathbf{k}=\lambda} t^{i+j+\lambda+1}(\ln t)^{\mu+1}\binom{\rho}{\mu}L_{j+1,\rho-\mu}(X_{\mathbf{r}\mathbf{l}})^{\mathbf{k}}B_{ij\mathbf{k}} +$$

$$+\sum_{i\geq 0,\ j\geq 0,\ \mathbf{k}\geq 0}\sum_{\rho\geq\mu}\sum_{\lambda=0}^{\infty}\sum_{\mathbf{r}:\ \mathbf{r}\bullet\mathbf{k}=\rho}\sum_{\mathbf{l}:\ \mathbf{l}\bullet\mathbf{k}=\lambda} t^{i+j+\lambda+1}(\ln t)^{\mu}\binom{\rho}{\mu}M_{j+1,\rho-\mu}(X_{\mathbf{r}\mathbf{l}})^{\mathbf{k}}B_{ij\mathbf{k}}$$

--- (5.10)

After suitable manipulations, we obtain the following system of equations for the unknown coefficients in the expansion (5.4):

$$X_{\mu n} = \Xi_{\mu n} + \sum_{\rho\geq\mu}\sum_{\substack{(\mathbf{r},\mathbf{l},\mathbf{k}):\\ \mathbf{l}\bullet\mathbf{k}\leq n-1,\ \mathbf{r}\bullet\mathbf{k}=\rho}}\sum_{i=0}^{n-\mathbf{l}\bullet\mathbf{k}-1}\binom{\rho}{\mu}L_{n-\mathbf{l}\bullet\mathbf{k}-i,\rho-\mu}(X_{\mathbf{r}\mathbf{l}})^{\mathbf{k}}A_{i,n-\mathbf{l}\bullet\mathbf{k}-i-1,\mathbf{k}} +$$

$$+\sum_{\rho\geq\mu}\sum_{\substack{(\mathbf{r},\mathbf{l},\mathbf{k}):\\ \mathbf{l}\bullet\mathbf{k}\leq n-1,\ \mathbf{r}\bullet\mathbf{k}=\rho}}\sum_{i=0}^{n-\mathbf{l}\bullet\mathbf{k}-1}\binom{\rho}{\mu}M_{n-\mathbf{l}\bullet\mathbf{k}-i,\rho-\mu}(X_{\mathbf{r}\mathbf{l}})^{\mathbf{k}}B_{i,n-\mathbf{l}\bullet\mathbf{k}-i-1,\mathbf{k}} +$$

$$+\chi(\mu\geq 1)\sum_{\rho\geq\mu-1}\sum_{\substack{(\mathbf{r},\mathbf{l},\mathbf{k}):\\ \mathbf{l}\bullet\mathbf{k}\leq n-1,\ \mathbf{r}\bullet\mathbf{k}=\rho}}\sum_{i=0}^{n-\mathbf{l}\bullet\mathbf{k}-1}\binom{\rho}{\mu-1}L_{n-\mathbf{l}\bullet\mathbf{k}-i,\rho-\mu+1}(X_{\mathbf{r}\mathbf{l}})^{\mathbf{k}}B_{i,n-\mathbf{l}\bullet\mathbf{k}-i-1,\mathbf{k}}$$

--- (5.11)

It is seen from (5.11) above that, in the case of logarithmic singularities, we do not obtain a recursive relation for the unknown coefficients but rather an infinite coupled polynomial system. The terms $(X_{\mathbf{r}\mathbf{l}})^{\mathbf{k}}$ depend only on those $X_{r_j l_j(j)}$ for which $k_j \geq 1$, since $(X_{r_j l_j(j)})^{k_j}$ becomes equal to 1 when $k_j = 0$; now, the condition $\mathbf{l}\bullet\mathbf{k} \leq n-1$ implies that, for those j for which $k_j \geq 1$, we have $l_j \leq n-1$, and therefore the right-hand side of (5.9) depends only on those $X_{r_j l_j(j)}$ that have $l_j \leq n-1$. Thus (5.11) is a sort of "partial" recursion, i.e. it is recursive with respect to n (out of the two subscripts in $X_{\mu n}$) but it is not recursive with respect to $\mu$ since the summation on the right-hand side is taken over $\rho \geq \mu$.

A similar situation occurs for Volterra systems with singularities of the type $(\ln s - \ln t)$. When f(t, s, x) has the form

$$f(t,s,x) = a(t,s,x) + (\ln s - \ln t)b(t,s,x), \text{ for } t > s;$$

$$a(t,s,x) = \sum_{i=0}^{\infty} \sum_{j=0}^{\infty} \sum_{\mathbf{k} \geq \mathbf{0}} t^i s^j x^{\mathbf{k}} A_{ij\mathbf{k}};$$

$$b(t,s,x) = \sum_{i=0}^{\infty} \sum_{j=0}^{\infty} \sum_{\mathbf{k} \geq \mathbf{0}} t^i s^j x^{\mathbf{k}} B_{ij\mathbf{k}}$$

--- (5.12)

with $\xi(t)$ given by (5.3), and we seek a solution x(t) in the form (5.4). After a series of calculations, of the same general nature as those in the case of singularities of the type $\ln(t-s)$, we obtain the following system:

$$X_{\mu n} = \Xi_{\mu n} + \sum_{\rho \geq \mu} \sum_{\substack{(\mathbf{r},\mathbf{l},\mathbf{k}): \\ \mathbf{l}\bullet\mathbf{k} \leq n-1, \, \mathbf{r}\bullet\mathbf{k}=\rho}} \sum_{i=0}^{n-\mathbf{l}\bullet\mathbf{k}-1} \binom{\rho}{\mu} L_{n-\mathbf{l}\bullet\mathbf{k}-i,\rho-\mu}(X_{\mathbf{r}\mathbf{l}})^{\mathbf{k}} A_{i,n-\mathbf{l}\bullet\mathbf{k}-i-1,\mathbf{k}} +$$

$$+ \sum_{\rho \geq \mu} \sum_{\substack{(\mathbf{r},\mathbf{l},\mathbf{k}): \\ \mathbf{l}\bullet\mathbf{k} \leq n-1, \, \mathbf{r}\bullet\mathbf{k}=\rho}} \sum_{i=0}^{n-\mathbf{l}\bullet\mathbf{k}-1} \binom{\rho}{\mu} L_{n-\mathbf{l}\bullet\mathbf{k}-i,\rho-\mu+1}(X_{\mathbf{r}\mathbf{l}})^{\mathbf{k}} B_{i,n-\mathbf{l}\bullet\mathbf{k}-i-1,\mathbf{k}}$$

--- (5.13)

Further, we note that, with related methods, we can obtain formal series solutions for Volterra equations with kernels having mixed singularities, of the form

$$f(t,s,x) = \sum_p \sum_q \sum_r c_{pqr}(t,s,x)(\ln(t-s))^p (\ln s - \ln t)^q (t-s)^{-\alpha_r}; \quad 0 < \alpha_r < 1;$$

$$c_{pqr}(t,s,x) = \sum_{i,j,\mathbf{k}} C_{pqrij\mathbf{k}} t^i s^j x^{\mathbf{k}}$$

--- (5.14)

A solution x(t) of a Volterra system with kernel given by (5.14) would then be sought in the form

$$x(t) = \sum_{p,r,\mu} \sum_{i_r: i_\mu > \mu\alpha_r - 1} (\ln t)^p t^{i_r - \mu\alpha_r} X_{pr\mu i_\mu}$$

--- (5.15)

Further, we can use the same ideas, but of course different technicalities, Volterra systems with kernels of the type

$$f(t,s,x) = \sum_{\lambda,\mu,\nu,\rho} c_{\lambda\mu\nu\rho}(t,s,x)(g_\lambda(t) - g_\lambda(s))^{p_\lambda} (h_\mu(t) - h_\mu(s))^{-\alpha_\mu} \cdot$$
$$\cdot (\varphi_\nu(t-s))^{q_\nu} (\psi_\rho(t-s))^{-\beta_\rho}$$

--- (5.16)

with functions $g_\lambda$, $h_\mu$, $\varphi_\nu$, $\psi_\rho$ satisfying certain technical conditions, where the exponents $p_\lambda$, $q_\nu$ are nonnegative integers, and the exponents $-\alpha_\mu$, $-\beta_\rho$ take values in (-1, 0). The functions $c_{\lambda\mu\nu\rho}(t,s,x)$ are analytic functions in the variables t, s, x.

as well as compositions of such functions, for example terms of the type $(\psi_\rho(g_\lambda(t) - g_\lambda(s)))^{-\gamma_{\rho\lambda}}$.

For the purposes of the present paper, we omit the formulation of the equations that must be satisfied by the coefficients $X_{pr\mu i_\mu}$ in these general cases.

## 6. Linear Volterra systems of the first kind.

In this section, we give sufficient conditions for a first-kind linear Volterra system to admit a Taylor series solution with recursive calculation of the coefficients.
We consider the system

$$\xi(t) + \int_0^t k(t,s)x(s)\,ds = 0$$

--- (6.1)

where

$$k(t,s) = \sum_{i=0}^{\infty} \sum_{j=0}^{\infty} t^i s^j K_{ij}, \quad \xi(t) = \sum_{i=0}^{\infty} t^i \Xi_i$$

--- (6.2)

where the series are assumed to have strictly positive radius of convergence.
It is plain that a necessary condition for existence of solutions of (6.1) is $\Xi_0 = 0$, since $\Xi_0 = \xi(0)$ and, for t=0, (6.1) gives $\xi(0) = 0$.
If we seek a solution

$$x(t) = \sum_{i=0}^{\infty} t^i X_i$$

--- (6.3)

then, as in section 3, we obtain

$$\Xi_n + \sum_{l=0}^{n-1} \sum_{i=0}^{n-l-1} \frac{1}{n-i} K_{i,n-i-l-1} X_l = 0, \text{ for } n \geq 1$$

--- (6.4)

If the matrix $K_{00}$ is nonsingular, then (6.4) is tantamount to a recursive formula for the coefficients $X_n$:

$$X_n = -(n+1)K_{00}^{-1}\left[\Xi_{n+1} + \sum_{l=0}^{n-1}\sum_{i=0}^{n-l} \frac{1}{n-i+1} K_{i,n-l-1}X_l\right]$$

--- (6.5)

More generally, we set

$$\mathbf{K}_{nj} := \sum_{i=0}^{j} \frac{1}{n-i+1} K_{i,j-i}, \text{ for } 0 \leq j \leq n$$

--- (6.6)

Then (6.4) becomes

$$\Xi_{n+1} + \sum_{j=0}^{n} \mathbf{K}_{nj} X_{n-j} = 0$$

--- (6.7)

Then it is seen from (6.7) that a possible scenario for the existence of recursive relations for the coefficients $X_i$ is the fulfillment of the following condition:

either $\mathbf{K}_{n0}$ is nonsingular for all n, or
$\mathbf{K}_{nj} = 0$ for $0 \leq j \leq j_0$ and $\mathbf{K}_{nj_0}$ is nonsingular for all n

--- (6.8)

In the case of the second possibility in (6.8), eq. (6.7) becomes

$$\Xi_{n+1} + \sum_{j=j_0}^{n} \mathbf{K}_{nj} X_{n-j} = 0$$

--- (6.9)

Clearly, (6.9) is possible only if $\Xi_{n+1} = 0$ for $n \leq j_0 - 1$, and then, for $n \geq j_0$, we obtain the recursive relationship

$$X_{n-j_0} = -\mathbf{K}_{nj_0}^{-1}\left[\Xi_{n+1} + \sum_{j=j_0+1}^{n} \mathbf{K}_{nj}X_{n-j}\right]$$

--- (6.10)